\documentclass[a4paper,11pt]{amsart}
\usepackage{amsmath,amssymb,amsfonts,latexsym}

\newtheorem{theorem}{Theorem}[section]

\newtheorem{corollary}[theorem]{Corollary}

\input{tcilatex}

\begin{document}
\title[\textbf{Dirichlet's type of twisted Eulerian polynomials}]{\textbf{%
Dirichlet's type of twisted Eulerian polynomials in connection with twisted
Eulerian-}$\mathcal{L}$\textbf{-function }}
\author[\textbf{S. Araci}]{\textbf{Serkan Araci}}
\address{\textbf{University of Gaziantep, Faculty of Science and Arts,
Department of Mathematics, 27310 Gaziantep, TURKEY}}
\email{\textbf{mtsrkn@hotmail.com}}
\author[\textbf{M. Acikgoz}]{\textbf{Mehmet Acikgoz}}
\address{\textbf{University of Gaziantep, Faculty of Science and Arts,
Department of Mathematics, 27310 Gaziantep, TURKEY}}
\email{\textbf{acikgoz@gantep.edu.tr}}

\begin{abstract}
In the present paper, we effect Dirichlet's type of twisted Eulerian
polynomials by using $p$-adic fermionic $q$-integral on the $p$-adic integer
ring. Also, we introduce some new interesting identities for them. As a
result of them, by using contour integral on the generating function of
Dirichlet's type of twisted Eulerian polynomials and so we define twisted
Eulerian-$\mathcal{L}$-function which interpolates of Dirichlet's type of
Eulerian polynomials at negative integers which we state in this paper.

\vspace{2mm}\noindent \textsc{2010 Mathematics Subject Classification.}
Primary 05A10, 11B65; Secondary 11B68, 11B73.

\vspace{2mm}

\end{abstract}

\maketitle

%%%%%%%%%%%%%%%%%%%%%%%%%%%%%%%%%%%%%%%%%%%%%%%%%%%%%%%%%%%%%%%%%%%

%%%%%%%%%%%%%%%%%%%%%%%%%%%%%%%%%%%%%%%%%%%%%%%%%%%%%%%%%%%%%%%%%%%

%%%%%%%%%%%%%%%%%%%%%%%%%%%%%%%%%%%%%%%%%%%%%%%%%%%%%%%%%%%%%%%%%%%

\section{\textbf{Introduction}}

%%%%%%%%%%%%%%%%%%%%%%%%%%%%%%%%%%%%%%%%%%%%%%%%%%%%%%%%%%%%%%%%%%%

As it is well-known, Eulerian polynomials, $\mathcal{A}_{n}\left( x\right) $%
, are given by means of the following exponential generating function:%
\begin{equation}
e^{\mathcal{A}\left( x\right) t}=\sum_{n=0}^{\infty }\mathcal{A}_{n}\left(
x\right) \frac{t^{n}}{n!}=\frac{1-x}{e^{t\left( 1-x\right) }-x}
\label{equation 19}
\end{equation}

where $\mathcal{A}^{n}\left( x\right) :=\mathcal{A}_{n}\left( x\right) $,
symbolically. Eulerian polynomials can find via the following recurrence
relation:%
\begin{equation}
\left( \mathcal{A}\left( t\right) +\left( t-1\right) \right) ^{n}-t\mathcal{A%
}_{n}\left( t\right) =\left\{ 
\begin{array}{cc}
1-t & \text{if }n=0 \\ 
0 & \text{if }n\neq 0\text{,}%
\end{array}%
\right.  \label{equation 20}
\end{equation}

(for details, see \cite{KIM}, \cite{Foata} and \cite{Araci 2}).

Imagine that $p$ be a fixed odd prime number. Throughout this paper, we use
the following notations. By $%
%TCIMACRO{\U{2124} }%
%BeginExpansion
\mathbb{Z}
%EndExpansion
_{p}$, we denote the ring of $p$-adic rational integers, $%
%TCIMACRO{\U{211a} }%
%BeginExpansion
\mathbb{Q}
%EndExpansion
$ denotes the field of rational numbers, $%
%TCIMACRO{\U{211a} }%
%BeginExpansion
\mathbb{Q}
%EndExpansion
_{p}$ denotes the field of $p$-adic rational numbers, and $%
%TCIMACRO{\U{2102} }%
%BeginExpansion
\mathbb{C}
%EndExpansion
_{p}$ denotes the completion of algebraic closure of $%
%TCIMACRO{\U{211a} }%
%BeginExpansion
\mathbb{Q}
%EndExpansion
_{p}$. Let $%
%TCIMACRO{\U{2115} }%
%BeginExpansion
\mathbb{N}
%EndExpansion
$ be the set of natural numbers and $%
%TCIMACRO{\U{2115} }%
%BeginExpansion
\mathbb{N}
%EndExpansion
^{\ast }=%
%TCIMACRO{\U{2115} }%
%BeginExpansion
\mathbb{N}
%EndExpansion
\cup \left\{ 0\right\} $.

The $p$-adic absolute value is defined by 
\begin{equation*}
\left\vert p\right\vert _{p}=\frac{1}{p}\text{.}
\end{equation*}%
Let $q$ be an indeterminate with $\left\vert q-1\right\vert _{p}<1$. Thus,
we give definition of $q$-integer of $x$, which is defined by 
\begin{equation*}
\left[ x\right] _{q}=\frac{1-q^{x}}{1-q}\text{ and }\left[ x\right] _{-q}=%
\frac{1-\left( -q\right) ^{x}}{1+q}\text{,}
\end{equation*}%
where we note that $\lim_{q\rightarrow 1}\left[ x\right] _{q}=x$ (see
[1-26]).

Let $UD\left( 
%TCIMACRO{\U{2124} }%
%BeginExpansion
\mathbb{Z}
%EndExpansion
_{p}\right) $ be the space of uniformly differentiable functions on $%
%TCIMACRO{\U{2124} }%
%BeginExpansion
\mathbb{Z}
%EndExpansion
_{p}$. For a positive integer $d$ with $\left( d,p\right) =1$, let 
\begin{equation*}
X=X_{d}=\lim_{\overleftarrow{n}}%
%TCIMACRO{\U{2124} }%
%BeginExpansion
\mathbb{Z}
%EndExpansion
/dp^{n}%
%TCIMACRO{\U{2124} }%
%BeginExpansion
\mathbb{Z}
%EndExpansion
=\underset{a=0}{\overset{dp-1}{\cup }}\left( a+dp%
%TCIMACRO{\U{2124} }%
%BeginExpansion
\mathbb{Z}
%EndExpansion
_{p}\right)
\end{equation*}

with%
\begin{equation*}
a+dp^{n}%
%TCIMACRO{\U{2124} }%
%BeginExpansion
\mathbb{Z}
%EndExpansion
_{p}=\left\{ x\in X\mid x\equiv a\left( \func{mod}dp^{n}\right) \right\}
\end{equation*}

where $a\in 
%TCIMACRO{\U{2124} }%
%BeginExpansion
\mathbb{Z}
%EndExpansion
$ satisfies the condition $0\leq a<dp^{n}$ and let $\mathcal{\sigma }%
:X\rightarrow 
%TCIMACRO{\U{2124} }%
%BeginExpansion
\mathbb{Z}
%EndExpansion
_{p}$ be the transformation introduced by the inverse limit of the natural
transformation%
\begin{equation*}
%TCIMACRO{\U{2124} }%
%BeginExpansion
\mathbb{Z}
%EndExpansion
/dp^{n}%
%TCIMACRO{\U{2124} }%
%BeginExpansion
\mathbb{Z}
%EndExpansion
\mapsto 
%TCIMACRO{\U{2124} }%
%BeginExpansion
\mathbb{Z}
%EndExpansion
/p^{n}%
%TCIMACRO{\U{2124} }%
%BeginExpansion
\mathbb{Z}
%EndExpansion
\text{.}
\end{equation*}

If $f$ is a function on $%
%TCIMACRO{\U{2124} }%
%BeginExpansion
\mathbb{Z}
%EndExpansion
_{p}$, then we will utilize the same notation to indicate the function $%
f\circ \mathcal{\sigma }$.

For a continuous function $f:X\rightarrow 
%TCIMACRO{\U{2102} }%
%BeginExpansion
\mathbb{C}
%EndExpansion
_{p},$ the $p$-adic $q$-integral on $%
%TCIMACRO{\U{2124} }%
%BeginExpansion
\mathbb{Z}
%EndExpansion
_{p}$ is defined by T. Kim in \cite{KIM2} and \cite{KIM4}, as follows:

\begin{equation}
I_{q}\left( f\right) =\int_{X}f\left( \upsilon \right) d\mu _{q}\left(
\upsilon \right) =\int_{%
%TCIMACRO{\U{2124} }%
%BeginExpansion
\mathbb{Z}
%EndExpansion
_{p}}f\left( \upsilon \right) d\mu _{q}\left( \upsilon \right)
=\lim_{n\rightarrow \infty }\frac{1}{\left[ p^{n}\right] _{q}}\sum_{\upsilon
=0}^{p^{n}-1}q^{\upsilon }f\left( \upsilon \right) \text{.}
\label{equation 2}
\end{equation}%
The bosonic integral is considered as the bosonic limit $q\rightarrow 1,$ $%
I_{1}\left( f\right) =\lim_{q\rightarrow 1}I_{q}\left( f\right) $. In \cite%
{Kim 2}, \cite{Kim 3} and \cite{Kim 4}, similarly, the $p$-adic fermionic
integration on $%
%TCIMACRO{\U{2124} }%
%BeginExpansion
\mathbb{Z}
%EndExpansion
_{p}$ \ is given by Kim as follows:%
\begin{equation}
I_{-q}\left( f\right) =\lim_{q\rightarrow -q}I_{q}\left( f\right) =\int_{%
%TCIMACRO{\U{2124} }%
%BeginExpansion
\mathbb{Z}
%EndExpansion
_{p}}f\left( \upsilon \right) d\mu _{-q}\left( \upsilon \right) \text{.}
\label{equation 3}
\end{equation}

By (\ref{equation 3}), we have the following well-known integral equation:%
\begin{equation}
q^{n}I_{-q}\left( f_{n}\right) +\left( -1\right) ^{n-1}I_{-q}\left( f\right)
=\left[ 2\right] _{q}\sum_{l=0}^{n-1}\left( -1\right) ^{n-1-l}q^{l}f\left(
l\right)  \label{equation 4}
\end{equation}

here let $f_{n}\left( x\right) $ be a translation with $f_{n}\left( x\right)
:=f\left( x+n\right) $. By (\ref{equation 4}), we readily derive the
following%
\begin{equation}
qI_{-q}\left( f_{1}\right) +I_{-q}\left( f\right) =\left[ 2\right]
_{q}f\left( 0\right) \text{.}  \label{equation 7}
\end{equation}

Replacing $q$ by $q^{-1}$ in (\ref{equation 7}),\ we easily derive the
following%
\begin{equation}
I_{-q^{-1}}\left( f_{1}\right) +qI_{-q^{-1}}\left( f\right) =\left[ 2\right]
_{q}f\left( 0\right) \text{.}  \label{equation 8}
\end{equation}

Recently, Kim $et$ $al.$ \cite{KIM} is considered $f(x)=e^{-x\left(
1+q\right) t}$ in (\ref{equation 8}), then they gave Witt's formula of
Eulerian polynomials as follows:

For $n\in 
%TCIMACRO{\U{2115} }%
%BeginExpansion
\mathbb{N}
%EndExpansion
^{\ast }$,%
\begin{equation}
I_{-q^{-1}}\left( x^{n}\right) =\frac{\left( -1\right) ^{n}}{\left(
1+q\right) ^{n}}\mathcal{A}_{n}\left( -q\right) \text{.}  \label{equation 15}
\end{equation}

In previous paper \cite{Araci 2}, Araci, Acikgoz and Gao are introduced
generating function of Dirichlet's type of Eulerian polynomials as%
\begin{equation}
\sum_{n=0}^{\infty }\mathcal{A}_{n,\chi }\left( -q\right) \frac{t^{n}}{n!}=%
\left[ 2\right] _{q}\sum_{l=0}^{d-1}\left( -1\right) ^{l}q^{d-l+1}\chi
\left( l\right) \frac{e^{-l\left( 1+q\right) t}}{e^{-d\left( 1+q\right)
t}+q^{d}}\text{.}  \label{equation 27}
\end{equation}

Also, they gave Witt's formula for Dirichlet's type of Eulerian polynomials,
as follows:%
\begin{equation*}
I_{-q^{-1}}\left( \chi \left( x\right) x^{n}\right) =\frac{\left( -1\right)
^{n}}{\left( 1+q\right) ^{n}}\mathcal{A}_{n,\chi }\left( -q\right) \text{.}
\end{equation*}

In this paper, we will also consider Dirichlet's type of twisted Eulerian
polynomials. By applying Mellin transformation to generating function of
Dirichlet's type of twisted Eulerian polynomials, we will describe twisted
Eulerian-$\mathcal{L}$-function. Next, by utilizing from Cauchy-Residue
theorem, we will see that twisted Eulerian-$\mathcal{L}$-function is related
to Dirichlet's type of twisted Eulerian polynomials at negative integers.

\section{\textbf{On the Dirichlet's type of twisted Eulerian polynomials}}

By using (\ref{equation 4}), it is easy to see that 
\begin{equation}
I_{-q^{-1}}\left( f_{d}\right) +q^{d}I_{-q^{-1}}\left( f\right) =\left[ 2%
\right] _{q}\sum_{l=0}^{d-1}\left( -1\right) ^{l}q^{d-l+1}f\left( l\right) 
\text{.}  \label{equation 28}
\end{equation}

Let $C_{p^{n}}=\left\{ \zeta \mid \zeta ^{p^{n}}=1\right\} $ be the cyclic
group of order $p^{n}$, and let 
\begin{equation*}
\mathcal{T}_{p}=\lim_{n\rightarrow \infty }C_{p^{n}}=C_{p^{\infty }}=%
\underset{n\geq 0}{\cup }C_{p^{n}}\text{,}
\end{equation*}

we want to note that $\mathcal{T}_{p}$ is locally constant space.

Let $\chi $ be a Dirichlet's character of conductor $d\left( =\text{odd}%
\right) $ and $\zeta \in \mathcal{T}_{p}$, then, taking $f(x)=\zeta ^{x}\chi
\left( x\right) e^{-x\left( 1+q\right) t}$ in (\ref{equation 28}), then it
is equality to%
\begin{gather*}
I_{-q^{-1}}\left( \zeta ^{x+d}\chi \left( x+d\right) e^{-\left( x+d\right)
\left( 1+q\right) t}\right) +q^{d}I_{-q^{-1}}\left( \zeta ^{x}\chi \left(
x\right) e^{-x\left( 1+q\right) t}\right) \\
=\left[ 2\right] _{q}\sum_{l=0}^{d-1}\left( -1\right) ^{l}q^{d-l+1}\chi
\left( l\right) e^{-l\left( 1+q\right) t}\zeta ^{l}
\end{gather*}

From this, we easily see that%
\begin{equation}
I_{-q^{-1}}\left( \zeta ^{x}\chi \left( x\right) e^{-x\left( 1+q\right)
t}\right) =\left[ 2\right] _{q}\sum_{l=0}^{d-1}\left( -1\right)
^{l}q^{d-l+1}\zeta ^{l}\chi \left( l\right) \frac{e^{-l\left( 1+q\right) t}}{%
\zeta ^{d}e^{-d\left( 1+q\right) t}+q^{d}}\text{.}  \label{equation 10}
\end{equation}

Now, we give definition of generating function of Dirichlet's type of
twisted Eulerian polynomials as follows:

%TCIMACRO{%
%\TeXButton{definition}{\begin{definition}
%Let $\mathcal{G}_{q,\zeta }\left( t\mid \chi \right) =\sum_{n=0}^{\infty }\mathcal{A}_{n,\chi ,\zeta }\left( -q\right) \frac{t^{n}}{n!}$ and $\zeta
%\in T_{p}$. Then, we define twisted Dirichlet's type of Eulerian polynomials by
%means of the following generating function:\begin{equation*}
%\mathcal{G}_{q,\zeta }\left( t\mid \chi \right) =\left[ 2\right]
%_{q}\sum_{l=0}^{d-1}\left( -1\right) ^{l}q^{d-l+1}\zeta ^{l}\chi \left(
%l\right) \frac{e^{-l\left( 1+q\right) t}}{\zeta ^{d}e^{-d\left( 1+q\right)
%t}+q^{d}}\text{.}
%\end{equation*}
%\end{definition}}}%
%BeginExpansion
\begin{definition}
Let $\mathcal{G}_{q,\zeta }\left( t\mid \chi \right) =\sum_{n=0}^{\infty }\mathcal{A}_{n,\chi ,\zeta }\left( -q\right) \frac{t^{n}}{n!}$ and $\zeta
\in T_{p}$. Then, we define twisted Dirichlet's type of Eulerian polynomials by
means of the following generating function:\begin{equation*}
\mathcal{G}_{q,\zeta }\left( t\mid \chi \right) =\left[ 2\right]
_{q}\sum_{l=0}^{d-1}\left( -1\right) ^{l}q^{d-l+1}\zeta ^{l}\chi \left(
l\right) \frac{e^{-l\left( 1+q\right) t}}{\zeta ^{d}e^{-d\left( 1+q\right)
t}+q^{d}}\text{.}
\end{equation*}
\end{definition}%
%EndExpansion

By considering the above definition, we readily derive the following
corollary.

\begin{corollary}
For any $\zeta \in \mathcal{T}_{p}$ and $n\in 
%TCIMACRO{\U{2115} }%
%BeginExpansion
\mathbb{N}
%EndExpansion
^{\ast }$, then we have%
\begin{equation*}
\mathcal{A}_{n,\chi ,\zeta }\left( -q\right) =\text{\textit{coefficient of} }%
\frac{t^{n}}{n!}\text{ \textit{of} }\left[ 2\right] _{q}\sum_{l=0}^{d-1}%
\left( -1\right) ^{l}q^{d-l+1}\zeta ^{l}\chi \left( l\right) \frac{%
e^{-l\left( 1+q\right) t}}{\zeta ^{d}e^{-d\left( 1+q\right) t}+q^{d}}\text{.}
\end{equation*}
\end{corollary}

By using (\ref{equation 10}) and Definition $1$, we procure the following
theorem which is the Witt's formula for Dirichlet's type of twisted Eulerian
polynomials.

\begin{theorem}
The following equality%
\begin{equation}
I_{-q^{-1}}\left( \zeta ^{x}\chi \left( x\right) x^{n}\right) =\int_{%
%TCIMACRO{\U{2124} }%
%BeginExpansion
\mathbb{Z}
%EndExpansion
_{p}}\zeta ^{x}\chi \left( x\right) x^{n}d\mu _{-q^{-1}}\left( x\right) =%
\frac{\left( -1\right) ^{n}}{\left( 1+q\right) ^{n}}\mathcal{A}_{n,\chi
,\zeta }\left( -q\right)  \label{equation 14}
\end{equation}%
holds true.
\end{theorem}

By using Definition 1, becomes%
\begin{eqnarray*}
\sum_{n=0}^{\infty }\mathcal{A}_{n,\chi ,\zeta }\left( -q\right) \frac{t^{n}%
}{n!} &=&\left[ 2\right] _{q}\sum_{l=0}^{d-1}\left( -1\right)
^{l}q^{d-l+1}\zeta ^{l}\chi \left( l\right) \frac{e^{-l\left( 1+q\right) t}}{%
\zeta ^{d}e^{-d\left( 1+q\right) t}+q^{d}} \\
&=&\left[ 2\right] _{q}\sum_{l=0}^{d-1}\left( -1\right) ^{l}q^{-l+1}\zeta
^{l}\chi \left( l\right) e^{-l\left( 1+q\right) t}\sum_{m=0}^{\infty }\left(
-1\right) ^{m}\zeta ^{md}q^{-md}e^{-md\left( 1+q\right) t} \\
&=&q\left[ 2\right] _{q}\sum_{m=0}^{\infty }\sum_{l=0}^{d-1}\left( -1\right)
^{l+md}\chi \left( l+md\right) q^{-\left( l+md\right) }\zeta
^{l+md}e^{-\left( l+md\right) \left( 1+q\right) t} \\
&=&q\left[ 2\right] _{q}\sum_{m=0}^{\infty }\left( -1\right) ^{m}\zeta
^{m}\chi \left( m\right) q^{-m}e^{-m\left( 1+q\right) t}\text{.}
\end{eqnarray*}

Thus, we get the following theorem.

\begin{theorem}
For $\zeta \in T_{p}$, then we have%
\begin{equation}
\mathcal{G}_{q,\zeta }\left( t\mid \chi \right) =\sum_{n=0}^{\infty }%
\mathcal{A}_{n,\chi ,\zeta }\left( -q\right) \frac{t^{n}}{n!}=\left[ 2\right]
_{q}\sum_{m=0}^{\infty }\frac{\left( -1\right) ^{m}\zeta ^{m}\chi \left(
m\right) e^{-m\left( 1+q\right) t}}{q^{m-1}}\text{.}  \label{equation 12}
\end{equation}
\end{theorem}

By using Taylor expansion of $e^{-m\left( 1+q\right) t}$ in (\ref{equation
12}), we discover the following.

\begin{theorem}
For $\zeta \in T_{p}$, then we have%
\begin{equation}
\frac{\left( -1\right) ^{n}}{q\left( 1+q\right) ^{n+1}}\mathcal{A}_{n,\chi
,\zeta }\left( -q\right) =\sum_{m=1}^{\infty }\frac{\left( -1\right)
^{m}\zeta ^{m}\chi \left( m\right) m^{n}}{q^{m}}\text{.}  \label{equation 13}
\end{equation}
\end{theorem}

From expressions of (\ref{equation 14}) and (\ref{equation 13}), we easily
effect the following corollary:

\begin{corollary}
For $\zeta \in T_{p}$ and $n\in 
%TCIMACRO{\U{2115} }%
%BeginExpansion
\mathbb{N}
%EndExpansion
$, then we procure the following%
\begin{equation*}
\lim_{n\rightarrow \infty }\sum_{x=0}^{p^{n}-1}\frac{\left( -1\right)
^{x}\zeta ^{x}\chi \left( x\right) x^{n}}{q^{x}}=2q^{2}\sum_{m=1}^{\infty }%
\frac{\left( -1\right) ^{m}\zeta ^{m}\chi \left( m\right) m^{n}}{q^{m}}\text{%
.}
\end{equation*}
\end{corollary}

Now, we want to show multiplication theorem for Dirichlet's type of twisted
Eulerian polynomials by using $p$-adic fermionic $q$-integral on $%
%TCIMACRO{\U{2124} }%
%BeginExpansion
\mathbb{Z}
%EndExpansion
_{p}$, as follows:

For each $\zeta \in \mathcal{T}_{p}$ and $n\in 
%TCIMACRO{\U{2115} }%
%BeginExpansion
\mathbb{N}
%EndExpansion
^{\ast }$, so we have 
\begin{align*}
& I_{-q^{-1}}\left( \zeta ^{x}\chi \left( x\right) x^{n}\right) \\
& =\int_{%
%TCIMACRO{\U{2124} }%
%BeginExpansion
\mathbb{Z}
%EndExpansion
_{p}}\zeta ^{x}\chi \left( x\right) x^{n}d\mu _{-q^{-1}}\left( x\right) \\
& =\lim_{m\rightarrow \infty }\frac{1}{\left[ dp^{m}\right] _{-q^{-1}}}%
\sum_{x=0}^{dp^{m}-1}\left( -1\right) ^{x}\zeta ^{x}\chi \left( x\right)
x^{n}q^{-x} \\
& =\frac{d^{n}}{\left[ d\right] _{-q^{-1}}}\sum_{a=0}^{d-1}\left( -1\right)
^{a}\zeta ^{a}\chi \left( a\right) q^{-a}\left( \lim_{m\rightarrow \infty }%
\frac{1}{\left[ p^{m}\right] _{-q^{-d}}}\sum_{x=0}^{p^{m}-1}\left( -1\right)
^{x}\zeta ^{dx}\left( \frac{a}{d}+x\right) ^{n}q^{-dx}\right) \\
& =\frac{d^{n}}{\left[ d\right] _{-q^{-1}}}\sum_{a=0}^{d-1}\left( -1\right)
^{a}\chi \left( a\right) \zeta ^{a}q^{-a}\int_{%
%TCIMACRO{\U{2124} }%
%BeginExpansion
\mathbb{Z}
%EndExpansion
_{p}}\left( \frac{a}{d}+x\right) ^{n}\zeta ^{dx}d\mu _{-q^{-d}}\left(
x\right) .
\end{align*}

Then, we can express Dirichlet's type of twisted Eulerian polynomials in
terms of $p$-adic fermionic $q$-integral on $%
%TCIMACRO{\U{2124} }%
%BeginExpansion
\mathbb{Z}
%EndExpansion
_{p}$, as follows:

\begin{theorem}
The following%
\begin{equation}
\frac{\left( -1\right) ^{n}}{\left( 1+q\right) ^{n}}\mathcal{A}_{n,\chi
,\zeta }\left( -q\right) =\frac{d^{n}}{\left[ d\right] _{-q^{-1}}}%
\sum_{a=0}^{d-1}\left( -1\right) ^{a}\chi \left( a\right) \zeta
^{a}q^{-a}\int_{%
%TCIMACRO{\U{2124} }%
%BeginExpansion
\mathbb{Z}
%EndExpansion
_{p}}\left( \frac{a}{d}+x\right) ^{n}\zeta ^{dx}d\mu _{-q^{-d}}\left(
x\right)  \label{equation 16}
\end{equation}%
is true.
\end{theorem}

Equation (\ref{equation 16}) seems to be interesting for evaluating at $q=1$%
, so we see that%
\begin{equation}
\frac{\left( -1\right) ^{n}}{2^{n}}\mathcal{A}_{n,\chi ,\zeta }\left(
-1\right) =d^{n}\sum_{a=0}^{d-1}\left( -1\right) ^{a}\chi \left( a\right)
\zeta ^{a}\int_{%
%TCIMACRO{\U{2124} }%
%BeginExpansion
\mathbb{Z}
%EndExpansion
_{p}}\left( \frac{a}{d}+x\right) ^{n}\zeta ^{dx}d\mu _{-1}\left( x\right) 
\text{.}  \label{equation 23}
\end{equation}

Next, with the help of Rim and Kim's paper \cite{Rim}, we can write the
following: 
\begin{equation}
E_{n,\zeta }\left( x\right) =\int_{%
%TCIMACRO{\U{2124} }%
%BeginExpansion
\mathbb{Z}
%EndExpansion
_{p}}\zeta ^{y}\left( x+y\right) ^{n}d\mu _{-1}\left( y\right)
\label{equation 21}
\end{equation}

where $\zeta \in \mathcal{T}_{p}$ and taking $x=0$ in the above equation, we
have $E_{n,\zeta }\left( 0\right) :=E_{n,\zeta }$, \ which is defined via
the following generating function:%
\begin{equation}
\sum_{n=0}^{\infty }E_{n,\zeta }\frac{t^{n}}{n!}=\frac{2\sum_{l=0}^{d-1}%
\left( -1\right) ^{l}\zeta ^{l}e^{lt}}{\zeta ^{d}e^{dt}+1},\text{ }%
\left\vert t\right\vert <\frac{\pi }{d}\text{.}  \label{equation 22}
\end{equation}

By expressions of (\ref{equation 23}), (\ref{equation 21}) and (\ref%
{equation 22}), we state the following corollary.

\begin{corollary}
For any $\zeta \in \mathcal{T}_{p}$ and $n\in 
%TCIMACRO{\U{2115} }%
%BeginExpansion
\mathbb{N}
%EndExpansion
^{\ast }$, then we obtain%
\begin{equation*}
\mathcal{A}_{n,\chi ,\zeta }\left( -1\right) =\left( -2d\right)
^{n}\sum_{a=0}^{d-1}\left( -1\right) ^{a}\chi \left( a\right) \zeta
^{a}E_{n,\zeta ^{d}}\left( \frac{a}{d}\right) \text{.}
\end{equation*}
\end{corollary}

\section{\textbf{On the twisted Eulerian-}$\mathcal{L}$\textbf{-function in }%
$\boldsymbol{%
%TCIMACRO{\U{2102} }%
%BeginExpansion
\mathbb{C}
%EndExpansion
}$}

The objective of this part is to describe twisted Eulerian-$\mathcal{L}$%
-function by applying Mellin transformation to the generating function of
Dirichlet's type of twisted Eulerian polynomials, which seem to be
interpolation function of Dirichlet's type of twisted Eulerian polynomials
at negative integers by using Cauchy-Residue theorem. We want to determine
that these functions are very useful in Number Theory, Complex Analysis and
Mathematical Physics, p-adic analysis and other areas cf. \cite{Cetin}, \cite%
{KIM4}, \cite{KIM5}, \cite{KIM 6}, \cite{Kim 6}, \cite{Jang}, \cite{Jang 1}, 
\cite{Cangul}, \cite{Cangul 1}, \cite{Araci}, \cite{Simsek}, \cite{Simsek 1}
and \cite{Simsek 2}. Thus, by (\ref{equation 12}), for $s\in 
%TCIMACRO{\U{2102} }%
%BeginExpansion
\mathbb{C}
%EndExpansion
$, we consider the following integral:%
\begin{equation}
\mathcal{L}_{E,\zeta }\left( s,\chi \right) =\frac{\int_{0}^{\infty }t^{s-1}%
\mathcal{G}_{q,\zeta }\left( t\mid \chi \right) dt}{\int_{0}^{\infty
}t^{s-1}e^{-t}dt}\text{.}  \label{equation 24}
\end{equation}

By (\ref{equation 24}), we can derive the following:%
\begin{eqnarray*}
\mathcal{L}_{E,\zeta }\left( s,\chi \right) &=&q\left[ 2\right]
_{q}\sum_{m=0}^{\infty }\left( -1\right) ^{m}\chi \left( m\right) \zeta
^{m}q^{-m}\left\{ \frac{\int_{0}^{\infty }t^{s-1}e^{-m\left( 1+q\right) t}dt%
}{\int_{0}^{\infty }t^{s-1}e^{-t}dt}\right\} \\
&=&\frac{q}{\left( 1+q\right) ^{s-1}}\sum_{m=1}^{\infty }\frac{\left(
-1\right) ^{m}\chi \left( m\right) \zeta ^{m}}{q^{m}m^{s}}\text{.}
\end{eqnarray*}

On the other hand, we easily see that%
\begin{equation}
=\sum_{n=0}^{\infty }\frac{\mathcal{A}_{n,\chi ,\zeta }\left( -q\right) }{n!}%
\left( \frac{\int_{0}^{\infty }t^{s-n+1}dt}{\int_{0}^{\infty }t^{s-1}e^{-t}dt%
}\right) \text{.}  \label{equation 25}
\end{equation}

So, we give definition of twisted Eulerian $\mathcal{L}$-function as follows:

\begin{definition}
Let $\zeta \in \mathcal{T}_{p}$ and $s\in 
%TCIMACRO{\U{2102} }%
%BeginExpansion
\mathbb{C}
%EndExpansion
$, then we define the following%
\begin{equation*}
\mathcal{L}_{E,\zeta }\left( s,\chi \right) =\frac{q}{\left( 1+q\right)
^{s-1}}\sum_{m=1}^{\infty }\frac{\left( -1\right) ^{m}\chi \left( m\right)
\zeta ^{m}}{q^{m}m^{s}}\text{.}
\end{equation*}
\end{definition}

By using (\ref{equation 24}) and (\ref{equation 25}), then, we derive the
following which seems to be interesting and useful for improving of Complex
Analysis and Theory of Analytic Numbers.

\begin{theorem}
The following equality holds true:%
\begin{equation*}
\mathcal{L}_{E,\zeta }\left( -n,\chi \right) =\left( -1\right) ^{n}\mathcal{A%
}_{n,\chi ,\zeta }\left( -q\right) \text{.}
\end{equation*}
\end{theorem}

%
%%%%%%%%%%%%%%%%%%%%%%%%%%%%%%%%%%%%%%%%%%%%%%%%%%%%%%%%%%%%%%%%

%%%%%%%%%%%%%%%%%%%%%%%%%%%%%%%%%%%%%%%%%%%%%%%%%%%%%%%%%%%%%%%%%%%

\end{document}